\documentclass[12pt]{article}
\usepackage{mathrsfs}
\usepackage{amsmath}
\usepackage{amssymb}
\usepackage{amsthm}
\usepackage{amsfonts}
\usepackage{amscd}
\usepackage[mathscr]{eucal}
\newcommand{\Z} {{\mathbb  Z}}
\newcommand{\Q}{{\mathbb  Q}}

\newcommand{\C}{{\mathbb  C}}

\newcommand{\R} {{\mathbb R}}
\newcommand{\f}{{\mathfrak f}}
\newcommand{\g}{{\mathfrak g}}
\newcommand{\h}{{\mathfrak h}}
\begin{document}
\parindent  25pt
\baselineskip  10mm
\textwidth  15cm    \textheight  23cm
\evensidemargin -0.06cm
\oddsidemargin -0.01cm

\title{ {A note on $ L(1) $ of Hecke $ L-$series associated
to the elliptic curves with CM by $ \sqrt{-3} $ }}
\author{\mbox{}
{ Derong Qiu }
\thanks{ \quad E-mail:
derong@mail.cnu.edu.cn } \\
(School of Mathematical Sciences, Institute of Mathematics \\
and Interdisciplinary Science, Capital Normal University, \\
Beijing 100048, P.R.China )  }

\date{}
\maketitle
\parindent  24pt
\baselineskip  10mm
\parskip  0pt

\par   \vskip 0.4cm

{\bf Abstract} \ Consider elliptic curves $ E:\ y^{2} = x^{3} +
D^{3} $ defined over the quadratic field $\ \Q(\sqrt{-3}) $. Hecke $
L-$series attached to $ E $ are studied, formulae for their values
at $ s=1, $ and bound of 3-adic valuations of these values are
given. These results are complementary to those in [Q] and [QZ], and
are consistent with the predictions of the conjecture of Birch and
Swinnerton-Dyer.
\par  \vskip  0.2 cm

\par  \vskip  0.1 cm

\par     \vskip  0.3 cm

\hspace{-0.6cm}{\bf 1. Introduction and statement of main results}

\par \vskip 0.2 cm

This note is a complement of [Q] and [QZ]. Let $ \tau = ({-1 +
\sqrt{-3}})/{2} $ be a primitive cubic root of unity and $ O_{K} =
\Z [\tau ] $ the ring of integers of the imaginary quadratic field $
K = \Q (\sqrt{-3}). $ In this note, we consider the elliptic curves
$$ E = E_{D^{3}}: \ y^{2} = x^{3} + D^{3}, \ \text{with} \
D = \pi _{1} \cdots \pi _{n}, \eqno{(1.1)} $$ where $ \pi _{k}
\equiv 1\ (\hbox{mod}\ 12) \ (k = 1, \cdots , n) $ are distinct
prime elements in $ O_{K}. $ Obviously, $ E $ has complex
multiplication by $ O_{K}. $ Let $ S = \{\pi _{1} \cdots \pi _{n}
\}. $ For any subset $ T $ of $ \{1, \cdots , n \}, $ denote $ D_{T}
= \prod _{k \in T} \pi _{k}, \ \widehat{D} _{T} = D / D_{T} $ and
put $ D_{\emptyset } = 1 $ when $ T = \emptyset $ (empty set). Let $
\psi _{D_{T} ^{3}} $ be the Hecke character (i.e.,
Gr$\ddot{o}$ssencharacter) of $ K $ attached to the elliptic curve $
E_{D_{T}^{3}}: \ y^{2} = x^{3} + D_{T}^{3}, $ and let $
L_{S}(\overline{\psi }_{D_{T} ^{3}}, \ s) $ be the Hecke $ L-$series
of $ \overline{\psi } _{D_{T} ^{3}}$ (the complex conjugate of $
{\psi } _{D_{T} ^{3}} $) with the Euler factors omitted at all
primes in $ S $ (for the definition of such Hecke $L-$series
attached to an elliptic curve, see [Sil2]). We have the following
result about the special value of $ L_{S}(\overline{\psi }_{D_{T}
^{3}}, \ s) $ at $ s = 1. $
\par \vskip 0.2cm

{\bf Theorem 1.1} \ Let $ \ D = \pi _{1} \cdots \pi _{n}, $ where $
\pi _{k} \equiv 1 \ (\hbox{mod} \ 12) $ are distinct prime elements
of $ \Z [ \tau ] \ ( k = 1, \cdots , n). $ Then, for any factor $
D_{T} $ of $ D $ and the corresponding Hecke character $ \psi
_{D_{T}^{3}}, $ we have
\begin{align*}
&- \frac{D}{\omega } \left(\frac{2}{D_{T}} \right) _{2}
L_{S}(\overline{\psi }_{D_{T} ^{3}}, \ 1) \\
&= \frac{\sqrt{3}}{4} \sum _{c \in \mathcal{C}}
\left(\frac{c}{D_{T}} \right) _{2} \frac{1}{\wp
\left(\frac{\sqrt{-3} c \omega }{D} \right) + \sqrt[3]{2}} -
\frac{\sqrt[3]{4}}{4 \sqrt{3}} \sum _{c \in \mathcal{C}}
\left(\frac{c}{D_{T}} \right) _{2}, \quad \quad \quad \quad (1.2)
\end{align*}
where $(-\ )_{_2} $ is the quadratic residue symbol in $ K, \
\mathcal{C}$ is any complete set of representatives of the
relatively prime residue classes of $ O_{K} $ modulo $ D, \ \wp (z)
$ is the Weierstrass $ \wp -$function satisfying $ \wp ^{\prime
}(z)^2 = 4 \wp (z)^{3} - 1 $ with period lattice  $ L_{\omega } =
\omega O_{K} $ (corresponding to the elliptic curve $ y^{2} = x^{3}
- \frac{1}{4} $) and $ \omega = 3.059908 \cdots  $ is an absolute
constant.
\par \vskip 0.1cm

There is much literature studying the special values $ L(1) $
associated to the CM elliptic curves (see e.g., [BSD], [Z$1\sim3$],
[Q], [QZ]). In [Q], a similar result of $ L(1) $ was obtained for
the Hecke character attached to some special elliptic curves $ y^{2}
= x^{3} - 2^{4} 3^{3} D^{3}. $ Now for the elliptic curves (1.1)
above, to obtain an explicit formula of $ L(1) $ need to overcome
more difficulties, especially in calculating the key values of
Weierstrass zeta function $ \zeta (z, \ L_{\omega}), $ Weierstrass $
\wp -$function $ \wp (z) (= \wp (z, \ L_{\omega})) $ and its
derivative $ \wp ^{\prime }(z) $ (see the proof of Theorem 1.1 in
the following).
 \par \vskip 0.1cm

Let $ {\Q _{2}} $ be the completion of $ \Q $ at the $ 2-adic $
valuation, $ \overline{\Q } $ and $ \overline{\Q _{2}}$ be the
algebraic closures of ${\Q }$ and $ \Q _{p} $ respectively, and let
$ v_{2} $ be the normalized $ 2-adic $ additive valuation of $
\overline{\Q _{2}} $ (i.e., $ v_{2}(2) = 1). $ Fix an isomorphic
embedding $ \overline{\Q} \hookrightarrow \overline{\Q _{2}}. $
Then, via this embedding, $ v_{2}(\alpha ) $ is defined for any
algebraic number $ \alpha $ in $ \overline{\Q }. $ The value $
v_{2}(\alpha ) $ for $ \alpha \in \overline{\Q } $ depends on the
choice of the embedding $ \overline{\Q} \hookrightarrow \overline{\Q
_{2}}, $ but this does not affect our discussion in this paper.
\par \vskip 0.1 cm

By Corollary 22 of [CW], we know that $ L(\overline{\psi } _{D^{3}},
\ 1) / \omega $ is an algebraic number, i.e., $ L(\overline{\psi
}_{D^{3}}, \ 1) / \omega \in \overline{\Q}. $ For its $2-$adic
valuation, we have
\par \vskip 0.2cm

{\bf Theorem 1.2.} \ Let $ D = \pi _{1} \cdots \pi _{n}, $ where $
\pi _{k} \equiv 1 \ (\hbox{mod} \ 12) $ are distinct prime elements
of $ \Z[\tau ] \ ( k = 1, \cdots , n), $ and let $ \psi _{D^{^{3}}}
$ be the Hecke character of $ \Q (\sqrt {-3}) $ attached to the
elliptic curve $ E_{D^{^{3}}}: y^{2} = x^{3} + D^{3}. $ Then, for
the $2-$adic valuation of $ L(\overline{\psi } _{D^{3}}, \ 1) /
\omega $ we have
$$  v_{2} \left(L(\overline{\psi } _{D^{3}}, \ 1) / \omega
\right) \geq n - 1.  $$

\par \vskip 0.3 cm

\hspace{-0.6cm}{\bf 2. Proofs of Theorems}

\vskip 0.2 cm

{\bf  Proof of Theorem 1.1.} \ For the elliptic curve $
E_{D_{T}^{3}}: \ y^{2} = x^{3} + D_{T}^{3}, $ it has complex
multiplication by $ O_{K}. $ Since the class number of $ K $ is 1,
the period lattice of $ E_{D_{T}^{3}}$ should be $ L_{T} = \omega
_{T} O_{K} $ for some $ \omega _{T} \in \C^{\times }. $ Let $ \omega
_{T} = \alpha _{T} \omega , \ \alpha _{T} \in \C^{\times }. $ By
Tate's algorithm [T], it is easy to show that the conductor of $
E_{D_{T}^{3}} $ is $ N_{E_{D_{T}^{3}}} = 12 D_{T}^{2}, $ and the
conductor of $ \psi _{D_{T}^{3}} $ is $ \f_{\psi _{D_{T}^{3}}} = (2
\sqrt{-3}D_{T}). $ In Prop.A of [QZ] (for a general form, see
Prop.5.5 in [GS]), putting $ k = 1, \ \h=O_{K}, \ \g =(2 \sqrt{-3}
D), \ \rho = \frac{\omega _{T}}{2 \sqrt{-3} D}, \ \phi = \psi
_{D_{T} ^{3}}, $ then the ray class field of $ K $ modulo  $ \g $ is
$ K((E_{D_{T} ^{3}}) _{\g}) $ (see the Lemma 4.7 in [GS]), and then
$$ \frac{\overline{\rho }}{|\rho | ^{2s}} L_{\g}(\overline{\psi }_{D_{T}
^{3}}, s) = \sum _{b \in \mathbf{B}} H_{1} \left( \frac{\psi _{D_{T}
^{3}}(b )\omega _{T}}{2 \sqrt{-3} D}, 0, s, L_{T} \right) \quad
(Re(s)> {3}/{2}) $$ with $ \mathbf{B} = \{(6 c + D): \ c \in
\mathcal{C} \}, $ such that $ \{\sigma _{b} : \ b \in \mathbf{B} \}
= \text{Gal} \left(K((E_{D_{T} ^{3}}) _{\g}) / K \right) \cong
\left(O_{K}/(2 \sqrt{-3} D) \right) ^{\times } / O_{K}^{\times }
\quad (\text{via Artin map}), $ where $ \mathcal{C} $ is as in
Theorem 1.1, a set of representatives of $ (O_{K}/(D)) ^{\times }. $
Then
$$ \frac{\overline{\rho }}{|\rho | ^{2s}} L_{\g}(\overline{\psi } _{D_{T}
^{3}}, s) = \sum _{c \in \mathcal{C}} H_{1} \left( \frac{\psi
_{D_{T} ^{3}}(6 c + D) \omega _{T}}{2 \sqrt{-3} D}, 0, s, \omega
_{T} O_{K} \right) \quad (Re(s)> {3}/{2}). $$ Note that $ H_{1}(z,
0, 1, L) $ could be analytically continued by the Eisenstein $
E^{*}-$function (see [W]): $ H_{1}(z, 0, 1, L) = E_{0,1} ^{*}(z, L)
= E_{1}^{*}(z, L). $ Hence we get
$$ \frac{2 \sqrt{-3} D}{\alpha _{T} \omega }
L_{\g}(\overline{\psi }_{D_{T} ^{3}}, 1) = \sum _{c \in \mathcal{C}}
E_{1} ^{*} \left(\psi _{D_{T} ^{3}} (6 c + D) \frac{ \alpha _{T}
\omega }{2 \sqrt{-3} D}, \alpha _{T} \omega O_{K} \right).
\eqno(2.1) $$ Since $ D \equiv 1 (\hbox{mod} \ 12), $ we have $ 6 c
+ D \equiv 1 \ (\hbox{mod} \ 6) $ for any $ c \in \mathcal{C}. $ In
particular, $ \left(\frac{2}{6 c + D} \right)_{3} = 1 $ (see [IR],
P.119). So by definition (see[Sil2], p.178),
$$ \psi _{D_{T} ^{3}}(6 c + D) = \overline{\left(\frac{4 D_{T}^{3}}{6 c
+ D} \right)_{6}}(6 c + D) = \left(\frac{D_{T}}{6 c + D}
\right)_{2}(6 c + D). $$ Moreover, by the quadratic reciprocity law
in $ K $ (see [Le], pp.256$\sim$260), we have $$
\left(\frac{D_{T}}{6 c + D} \right)_{2} = \left(\frac{6 c +
D}{D_{T}} \right)_{2} = \left(\frac{6 c}{D_{T}} \right)_{2} =
\left(\frac{-2 \cdot (\sqrt{-3})^{2} c}{D_{T}} \right)_{2} =
\left(\frac{- 2 c}{D_{T}} \right)_{2} = \left(\frac{2 c}{D_{T}}
\right)_{2}, $$ the last equality holds because $ \left(\frac{-
1}{D_{T}} \right)_{2} = 1 $ (see [Le], p.111). Therefore, by (2.1)
above, and note that $ L_{\g}(\overline{\psi }_{D_{T} ^{3}}, 1) =
L_{S}(\overline{\psi }_{D_{T} ^{3}}, 1), $ we obtain
$$ \frac{2 \sqrt{-3} D}{\alpha _{T} \omega } L_{S}
(\overline{\psi }_{D_{T} ^{3}}, 1) = \sum _{c \in \mathcal{C}} E_{1}
^{*} \left(\left(- \frac{\sqrt{-3} c \omega }{D} - \frac{\sqrt{-3}
\omega }{6} \right) \alpha _{T} \left(\frac{2 c}{D_{T}} \right)
_{2}, \ \alpha _{T} \omega O_{K} \right). \eqno(2.2) $$ Let $
\lambda = - \alpha _{T} \left(\frac{2 c}{D_{T}} \right) _{2}, $ then
$ \alpha _{T} \omega O_{K} = \lambda \omega O_{K} = \lambda
L_{\omega }. $ By formula $ E_{1} ^{*}(\lambda z, \lambda L) =
\lambda ^{-1} E_{1} ^{*}(z, L), $ we obtain
\begin{align*} &E_{1}^{*} \left(\left(- \frac{\sqrt{-3} c \omega }{D} -
\frac{\sqrt{-3} \omega }{6} \right) \alpha _{T} \left(\frac{2
c}{D_{T}} \right) _{2}, \ \alpha _{T} \omega O_{K} \right) \\
&= E_{1} ^{*} \left(\left( \frac{\sqrt{-3} c \omega }{D} +
\frac{\sqrt{-3} \omega }{6} \right) \lambda , \ \lambda L_{\omega }
\right) = - \alpha _{T} ^{-1} \left( \frac{2 c}{D_{T}} \right)_{2}
E_{1} ^{*} \left(\frac{\sqrt{-3} c \omega }{D} + \frac{\sqrt{-3}
\omega }{6}, \ L_{\omega } \right).
\end{align*} So by (2.2) above, we get
$$ - \frac{D}{\omega } \left(\frac{2}{D_{T}} \right) _{2}
L_{S} (\overline{\psi }_{D_{T}^{3}}, 1) = \frac{1}{2 \sqrt{-3}} \sum
_{c \in \mathcal{C}} \left(\frac{c}{D_{T}} \right) _{2} E_{1}^{*}
\left(\frac{\sqrt{-3} c \omega }{D} + \frac{\sqrt{-3} \omega }{6}, \
L_{\omega } \right). \eqno(2.3) $$ By [QZ], it is easy to see that
$$ E_{1}^{*}(z, L_{\omega }) = \zeta (z, L_{\omega }) -
\frac{2 \pi \overline{z}}{\sqrt{3} \omega ^{2}}, \eqno(2.4) $$ where
$ \zeta (z, L) = \frac{1}{z} + \sum _{\alpha \in L-\{0 \}}
\left(\frac{1}{z - \alpha } + \frac{1}{\alpha } + \frac{z}{\alpha
^{2}} \right) $ is the Weierstrass Zeta-function, an odd function,
i.e., $ \zeta (- z, L) = - \zeta (z, L) $ (see [Sil 2]). By the
addition formula (see [Law])
$$ \zeta (z_{1} + z_{2}, L_{\omega }) = \zeta (z_{1}, L_{\omega }) +
\zeta (z_{2}, L_{\omega }) + \frac{1}{2} \frac{\wp ^{\prime }(z_{1})
- \wp ^{\prime }(z_{2})} {\wp (z_{1}) - \wp (z_{2})}, \quad \text{we
obtain} $$
\begin{align*}
&\zeta \left(\frac{\sqrt{-3} c \omega }{D} + \frac{\sqrt{-3} \omega
}{6}, \ L_{\omega } \right) \\
&= \zeta \left(\frac{\sqrt{-3} c \omega }{D},  L_{\omega } \right) +
\zeta \left(\frac{\sqrt{-3} \omega }{6}, L_{\omega } \right) +
\frac{1}{2} \frac{\wp ^{\prime } \left(\frac{\sqrt{-3} c \omega }{D}
\right) - \wp ^{\prime } \left(\frac{\sqrt{-3} \omega }{6}
\right)}{\wp \left(\frac{\sqrt{-3} c \omega }{D} \right) - \wp
\left(\frac{\sqrt{-3} \omega }{6} \right)}.  \quad \quad (2.5)
\end{align*}
Now we compute the values of $ \zeta \left(\frac{\sqrt{-3} \omega
}{6}, L_{\omega } \right), \ \wp \left(\frac{\sqrt{-3} \omega }{6}
\right) $ and $ \wp ^{\prime } \left(\frac{\sqrt{-3} \omega }{6}
\right). $ Note that $ \sqrt{-3} = 1 + 2 \tau , \ \tau O_{K} =
O_{K}, \ \tau L_{\omega } = L_{\omega }, \wp (\tau z, L_{\omega }) =
\tau \wp (z, L_{\omega }), \ \wp ^{\prime }(\tau z, L_{\omega }) =
\wp ^{\prime }(z, L_{\omega }) $ and $ \zeta (\tau z, L_{\omega }) =
\tau ^{2} \zeta (z, L_{\omega }) $ (see [La], p.16, p.240). Also by
[St] and [QZ], we know that
\begin{align*}
&\wp \left(\frac{\omega }{3}, L_{\omega } \right) = 1, \quad \wp
^{\prime } \left(\frac{\omega }{3}, L_{\omega } \right) = -
\sqrt{3}; \quad \wp ^{\prime \prime } \left(\frac{\omega }{3},
L_{\omega } \right) = 6, \quad \zeta \left(\frac{\omega }{2},
L_{\omega } \right) = \frac{\pi } {\sqrt{3} \omega }, \\
&\zeta \left(\frac{\omega }{3}, L_{\omega } \right) = \frac{2 \pi }
{3 \sqrt{3} \omega } + \frac{1}{\sqrt{3}}, \quad \zeta \left(\frac{2
\omega }{3}, L_{\omega } \right) = \frac{4 \pi } {3 \sqrt{3} \omega
} - \frac{1}{\sqrt{3}}. \quad (2.6)
\end{align*}
For $ O_{K} = \Z [\tau ], $ it is easy to see that the Eisenstein
series $ G_{2k}(O_{K}) $ is a real number for each positive integer
$ k \geq 2. $ So by the Laurent series expansion $ \wp (z, O_{K}) =
z^{-2} + \sum _{k = 1}^{\infty } (2 k + 1) G_{2k + 2}(O_{K}) z^{2k}
$ (see [Sil1], p.169), it is easy to see that $ \wp (\frac{1}{2},
O_{K}) \in \R, $ a real number, so $ \wp \left(\frac{\omega }{2},
L_{\omega } \right) = \omega ^{-2} \wp (\frac{1}{2}, O_{K}) \in \R.
$ Then, since $ (\wp (\frac{\omega }{2}, L_{\omega }), \ \frac{1}{2}
\wp ^{\prime }(\frac{\omega }{2}, L_{\omega })) $ is a point of
order $ 2 $ of the elliptic curve $ y^{2} = x^{3} - \frac{1}{4} $
mentioned above, one can easily obtain that $$ \wp ^{\prime }
\left(\frac{\omega }{2}, L_{\omega } \right) = 0, \ \wp
\left(\frac{2 \omega }{3}, L_{\omega } \right) = 1, \ \wp ^{\prime }
\left(\frac{2\omega }{3}, L_{\omega } \right) = \sqrt{3}, \ \wp
\left(\frac{\omega }{2}, L_{\omega } \right) =
\frac{\sqrt[3]{2}}{2}. \eqno(2.7) $$ So by the addition formula of $
\zeta (z, L_{\omega }) $ above, we get
\begin{align*} &\zeta \left(\frac{5 \omega}{6}, L_{\omega } \right)
= \zeta \left(\frac{\omega }{2} + \frac{\omega }{3}, L_{\omega }
\right) = \zeta \left(\frac{\omega }{2}, L_{\omega } \right) + \zeta
\left(\frac{\omega }{3}, L_{\omega } \right) + \frac{1}{2} \cdot
\frac{\wp ^{\prime }(\frac{\omega}{2}) - \wp ^{\prime
}(\frac{\omega}{3})} {\wp (\frac{\omega}{2}) - \wp (\frac{\omega}{3})} \\
&= \frac{5 \pi }{3 \sqrt{3} \omega } + \frac{1}{\sqrt{3}} +
\frac{\sqrt{3}}{\sqrt[3]{2} - 2}. \quad \quad \quad \quad (2.8)
\end{align*}
Moreover, for any $ \alpha \in L_{\omega }, $ we have
$$ \zeta (z + \alpha , L_{\omega }) -
\zeta (z, L_{\omega }) = \eta (\alpha , L_{\omega }) = \alpha
s_{2}(L_{\omega }) + \overline{\alpha } A(L_{\omega })^{-1} =
\frac{2 \pi \overline{\alpha }}{\sqrt{3} \omega ^{2}}
$$ because $ s_{2}(L_{\omega }) = \frac{2}{\omega}
\zeta (\frac{\omega}{2}, L_{\omega }) - \frac{2 \pi }{\sqrt{3}
\omega ^{2}} = 0 $ and $ A(L_{\omega }) = \frac{\sqrt{3} \omega
^{2}}{2 \pi } $ (see [QZ]). Putting $ z = - \frac{\omega }{6} $ and
$ \alpha = \omega , $ then we obtain
\begin{align*} &\zeta
\left(\frac{5 \omega}{6}, L_{\omega } \right) + \zeta
\left(\frac{\omega}{6}, L_{\omega } \right) = \frac{2 \pi }{\sqrt{3}
\omega }. \quad \text{So by (2.8), we get} \\
&\zeta \left(\frac{\omega}{6}, L_{\omega } \right) = \frac{\pi }{3
\sqrt{3} \omega } - \frac{1}{\sqrt{3}} - \frac{\sqrt{3}}{\sqrt[3]{2}
- 2}. \quad \quad (2.9)
\end{align*}
Also by taking $ u = \frac{2\omega }{3} $ and $ v = \frac{\omega
}{6} $ in the formula (see [Law], p.161)
\begin{align*} &\zeta (u + v, L_{\omega }) + \zeta (u - v, L_{\omega })
- 2 \zeta (u, L_{\omega }) = \frac{\wp ^{\prime }(u)}{\wp (u) - \wp
(v)}, \quad \quad (2.10) \\
&\text{we get} \quad \zeta \left(\frac{2 \omega}{3} + \frac{\omega
}{6}, L_{\omega } \right) + \zeta \left(\frac{2 \omega}{3} -
\frac{\omega }{6}, L_{\omega }\right) - 2 \zeta \left(\frac{2
\omega}{3}, L_{\omega } \right) = \frac{\wp ^{\prime }(\frac{2
\omega}{3})}{\wp (\frac{2 \omega}{3}) - \wp (\frac{\omega }{6})}, \\
&\text{which implies} \quad \wp \left(\frac{\omega }{6} \right) = 1
+ \sqrt[3]{2} + \sqrt[3]{4}. \quad \quad (2.11)
\end{align*}
Then by taking $ u = \frac{\omega }{6} $ and $ v = \frac{\omega }{3}
$ in the formula (2.10) above, we get
$$ \wp ^{\prime } \left(\frac{\omega}{6} \right) =
- \sqrt{3} \left(3 + 2 \cdot \sqrt[3]{2} + 2 \cdot \sqrt[3]{4}
\right). \eqno(2.12) $$ Now, by substituting these values into the
addition formula of $ \wp (z) $ (see [Law], p.162), we have
\begin{align*} &\wp \left(\frac{\sqrt{-3} \omega }{6} \right) = \wp
\left(\frac{\omega}{6} + \frac{\tau \omega }{3} \right) =
\frac{1}{4} \left( \frac{\wp ^{\prime } \left(\frac{\omega}{6}
\right) - \wp ^{\prime } \left(\frac{\tau \omega}{3} \right)}{\wp
\left(\frac{\omega}{6} \right) - \wp \left(\frac{\tau \omega}{3}
\right)} \right)^{2} - \wp \left(\frac{\omega }{6} \right) - \tau
\wp \left(\frac{\omega }{3} \right) = - \sqrt[3]{2},  \\
&\text{that is} \quad \wp \left(\frac{\sqrt{-3} \omega }{6} \right)
= - \sqrt[3]{2}; \quad \quad (2.13)
\end{align*}
\begin{align*} &\wp \left(\frac{\omega }{\sqrt{-3}} \right) = \wp
\left(\frac{\omega }{3} + \frac{2 \tau \omega }{3} \right) =
\frac{1}{4} \left( \frac{\wp ^{\prime } \left(\frac{\omega}{3}
\right) - \wp ^{\prime } \left(\frac{2 \tau \omega}{3} \right)}{\wp
\left(\frac{\omega}{3} \right) - \wp \left(\frac{2 \tau \omega}{3}
\right)} \right)^{2} - \wp \left(\frac{\omega }{3} \right) - \tau
\wp \left(\frac{2 \omega }{3} \right) = 0,  \\
&\text{that is} \quad \wp \left(\frac{\omega }{\sqrt{-3}} \right) =
0. \quad \quad (2.14)
\end{align*}
Next, by putting $ u = \frac{\omega }{6} $ and $ v = \frac{\tau
\omega }{3} $ into the following formula (see [Law], p.183, Exer.
15) $$ \frac{\wp ^{\prime }(u) - \wp ^{\prime }(v)}{\wp (u) - \wp
(v)} = \frac{\wp ^{\prime }(v) + \wp ^{\prime }(u + v)}{\wp (v) -
\wp (u + v)},  \ \ \text{we obtain} \ \ \wp ^{\prime }
\left(\frac{\sqrt{-3} \omega }{6} \right) = - 3 \sqrt{-1}.
\eqno(2.15) $$ Again by the addition formula of $ \zeta (z,
L_{\omega }) $ above, we get
\begin{align*} \zeta \left(\frac{\sqrt{-3} \omega }{6}, L_{\omega }
\right) &= \zeta \left(\frac{\omega }{6} + \frac{\tau \omega }{3},
L_{\omega } \right) = \zeta \left(\frac{\omega }{6}, L_{\omega }
\right) + \zeta \left(\frac{\tau \omega }{3}, L_{\omega } \right) +
\frac{1}{2} \cdot \frac{\wp ^{\prime }(\frac{\omega }{6}) - \wp
^{\prime }(\frac{\tau \omega }{3})}
{\wp (\frac{\omega}{6}) - \wp (\frac{\tau \omega }{3})} \\
&= - \frac{\pi \cdot \sqrt{-1}}{3 \omega } - \frac{\sqrt{-1}}{2}
\cdot \sqrt[3]{4}, \quad \text{and then by (2.5) we obtain}
\end{align*}
\begin{align*} &\zeta
\left(\frac{\sqrt{-3} c \omega }{D} + \frac{\sqrt{-3} \omega }{6},
L_{\omega } \right) \\
&= \zeta \left(\frac{\sqrt{-3} c \omega }{D}, L_{\omega } \right) -
\frac{\pi \cdot \sqrt{-1}}{3 \omega } - \frac{\sqrt{-1}}{2} \cdot
\sqrt[3]{4} + \frac{1}{2} \cdot \frac{\wp ^{\prime }(\frac{\sqrt{-3}
c \omega }{D}) + 3 \sqrt{-1}} {\wp (\frac{\sqrt{-3} c \omega }{D}) +
\sqrt[3]{2}}. \quad (2.16)
\end{align*}
Substituting it into (2.4) and (2.3), we obtain
\begin{align*}
&- \frac{D}{\omega } \left(\frac{2}{D_{T}} \right) _{2}
L_{S}(\overline{\psi }_{D_{T} ^{3}}, \ 1) = \frac{\sqrt{3}}{4} \sum
_{c \in \mathcal{C}} \left(\frac{c}{D_{T}} \right) _{2} \frac{1}{\wp
\left(\frac{\sqrt{-3} c \omega }{D} \right) + \sqrt[3]{2}} -
\frac{\sqrt[3]{4}}{4 \sqrt{3}} \sum _{c \in \mathcal{C}}
\left(\frac{c}{D_{T}} \right) _{2} \\
&+ \frac{1}{2 \sqrt{-3}} \sum _{c \in \mathcal{C}}
\left(\frac{c}{D_{T}} \right) _{2} \left( \zeta
\left(\frac{\sqrt{-3} c \omega }{D}, L_{\omega } \right) +
\frac{1}{2} \cdot \frac{\wp ^{\prime }(\frac{\sqrt{-3} c \omega
}{D})} {\wp (\frac{\sqrt{-3} c \omega }{D}) + \sqrt[3]{2}} + \frac{2
\pi \sqrt{-1}}{\omega } \cdot \frac{\overline{c}}{\overline{D}}
\right).
\end{align*}
Since $ D = \pi _{1} \cdots \pi _{n} $ with $ \pi _{k} \equiv 1 \
(\hbox{mod} \ 12), $ so we may choose the set  $ \mathcal{C} $ in
such a way that $ - c \in \mathcal{C} $ when $ c \in \mathcal{C}. $
Obviously  $ \left({-c}/{D_{T}} \right) _{2} = \left({c}/{D_{T}}
\right) _{2}. $ Also since  $ \zeta (z, L_{\omega }) $ and $ \wp
^{\prime }(z, L_{\omega }) $ are odd functions , and $ \wp (z,
L_{\omega })$ is an even function, so
$$ \sum _{c \in \mathcal{C}} \left(\frac{c}{D_{T}} \right) _{2}
\zeta \left(\frac{\sqrt{-3} c \omega }{D}, L_{\omega } \right) =
\sum _{c \in \mathcal{C}} \left(\frac{c}{D_{T}} \right) _{2}
\frac{\wp ^{\prime }(\frac{\sqrt{-3} c \omega }{D})} {\wp
(\frac{\sqrt{-3} c \omega }{D}) + \sqrt[3]{2}} = \sum _{c \in
\mathcal{C}} \left(\frac{c}{D_{T}} \right) _{2}
\frac{\overline{c}}{\overline{D}} = 0. $$ Therefore $$ -
\frac{D}{\omega } \left(\frac{2}{D_{T}} \right) _{2}
L_{S}(\overline{\psi }_{D_{T} ^{3}}, \ 1) = \frac{\sqrt{3}}{4} \sum
_{c \in \mathcal{C}} \left(\frac{c}{D_{T}} \right) _{2} \frac{1}{\wp
\left(\frac{\sqrt{-3} c \omega }{D} \right) + \sqrt[3]{2}} -
\frac{\sqrt[3]{4}}{4 \sqrt{3}} \sum _{c \in \mathcal{C}}
\left(\frac{c}{D_{T}} \right) _{2}. $$ This proves Theorem 1.1.
\quad $ \Box $
\par \vskip 0.2 cm

{\bf Remark 2.1.} (1) \ It follows from the above proof that Theorem
1.1 holds for all $ D = \pi _{1} \cdots \pi _{n} $ with $ \pi _{k}
\equiv 1 \ (\hbox{mod} \ 4 \sqrt{-3}). $ \\
(2) \ In particular, by taking $ D = 1 $ in the formula of Theorem
1.1, it is ease to see that $ L(E_{1} / \Q, 1) = L(\overline{\psi
}_{1}, 1) = \frac{\sqrt[3]{4}}{4 \sqrt{3}} \cdot \omega  $ for the
elliptic curve $ E_{1} : y^{2} = x^{3} + 1. $
\par \vskip 0.2 cm

{\bf Lemma 2.2.} \ For the Weierstrass $ \wp -$function $ \wp (z,
L_{\omega }) $ in Theorem 1.1 and any $ c \in \mathcal{C}, $ we have
$$ v_{2} \left(\wp \left(\frac{\sqrt{-3} c \omega }{D}, \
L_{\omega } \right) + \sqrt[3]{2} \right) = 0. $$
\par \vskip 0.4cm

{\bf Proof.} \ Taking $ r = 1, \ \gamma = 1, \ \Delta = D, \ \beta =
\sqrt{-3} c $ and $ \lambda = \frac{1}{2}(1 - 3^{1 - r}) = 0 $ in
the lemmas 2 and 1 in [St], by the above (2.14), it then follows
that $$ v_{2} \left(\wp \left( \frac{\sqrt{-3} c \omega }{D} \right)
\right) = 0, \quad \text{so} \quad v_{2} \left(\wp \left(
\frac{\sqrt{-3} c \omega }{D} \right) + \sqrt[3]{2} \right) = 0. $$
The proof is completed. \quad $ \Box $
\par \vskip 0.2cm

{\bf Proof of Theorem 1.2.} \ Add up the two sides of the formula in
Theorem 1.1 over all subsets $ T $ of $ \{1, \cdots , n \}, $ we
obtain
$$ - \sum _{T} \frac{D}{\omega } \left(\frac{2}{D_{T}} \right) _{2}
L_{S}(\overline{\psi } _{D_{T} ^{3}}, 1) = \frac{\sqrt{3}}{4} \sum
_{c \in \mathcal{C}} \frac{1}{\wp \left(\frac{\sqrt{-3} c \omega
}{D} \right) + \sqrt[3]{2}} \sum _{T} \left(\frac{c}{D_{T}} \right)
_{2} - \frac{\sqrt[3]{4}}{4 \sqrt{3}} \cdot \sharp \mathcal{C}.
\eqno (2.17) $$ By assumption,
$$ v_{2} \left(\frac{\sqrt[3]{4}}{4 \sqrt{3}} \cdot \sharp
\mathcal{C} \right) = v_{2} \left(\frac{\sqrt[3]{4}}{4 \sqrt{3}}
\cdot \prod \limits _{k = 1}^{n}(\pi _{k} \overline{\pi
_{k}}-1)\right) \geq \frac{2}{3} - 2 + 2 n = 2n - \frac{4}{3}. \quad
(n\geq 1) $$ Note that by our choice $ -c \in \mathcal{C} $ when $ c
\in \mathcal{C}, $ and $ \left(\frac{- c}{D_{T}} \right) _{2} =
\left(\frac{c}{D_{T}} \right) _{2}, $ so by Lemma 2.2 we know that
the first term in the right side of (2.17) has $ 2-$adic valuation $
\geq -2 + 1 + n = n -1. $ Therefore
$$ v_{2} \left(\sum _{T} \frac{D}{\omega }
\left(\frac{2}{D_{T}} \right) _{2} L_{S}(\overline{\psi } _{D_{T}
^{3}}, 1) \right) \geq n - 1. \eqno(2.18) $$ By definition, we know
that, if $ T = \{1, \cdots , n \}, $ then $ L_{S}(\overline{\psi
}_{D_{T} ^{3}}, 1) = L(\overline{\psi } _{D ^{3}}, 1); $ and if $ T
= \emptyset , $ then $ L_{S}(\overline{\psi }_{D_{T} ^{3}}, 1) =
L_{S}(\overline{\psi }_{1}, 1) = L(\overline{\psi } _{1}, 1) \prod
\limits _{k = 1}^{n} \left(1- \frac{1}{\pi _{k}} \right)=
\frac{\sqrt[3]{4}}{4 \sqrt{3}} \cdot \omega \prod \limits _{k =
1}^{n} \left(1 - \frac{1}{\pi _{k}} \right) $ (see the above Remark
2.1.(2)). So we have
$$ v_{2} \left(L_{S}(\overline{\psi } _{1}, 1)/\omega \right)
\geq -\frac{4}{3} + 2 n \geq n - 1 \quad (\hbox{Since} \ v_{2} (\pi
_{k} - 1) \geq 2). \eqno(2.19) $$ Now we use induction method on $n$
to prove $ v_{2} \left(L (\overline{\psi }_{D^{3}}, 1)/ \omega
\right)\geq n - 1. $ When $ n = 1, D = \pi _{1}, \ v_{2}
\left(L_{S}(\overline{\psi } _{1}, 1) /\omega \right) \geq
-\frac{4}{3} + 2 = \frac{2}{3}. $ Also by taking $ n = 1 $ in
(2.18),
$$ v_{2} \left(\frac{\pi _{1}}{\omega } \left(\frac{2}{D_{\emptyset
}} \right)_{2} L_{S}(\overline{\psi }_{1}, 1) + \frac{\pi
_{1}}{\omega } \left(\frac{2}{\pi _{1}} \right)_{2} L(\overline{\psi
}_{\pi _{1} ^{3}}, 1) \right) \geq 1 - 1 =0. $$ So $ v_{2}
\left(L(\overline{\psi }_{\pi _{1} ^{3}}, 1)/\omega \right)= v_{2}
\left(\frac{\pi _{1}}{\omega } \left(\frac{2}{\pi _{1}} \right) _{2}
L(\overline{\psi } _{\pi _{1} ^{3}}, 1) \right) \geq 0. $ \ Assume
our conclusion is true for $ 1, 2, \cdots , n-1, $ and consider the
case $ n, \ D = \pi _{1} \cdots \pi _{n}. $ For any non-trivial
subset $ T $ of $ \{1, \cdots , n \}, $ denote $ t = t(T) =\sharp T,
$ by definition, we have
\begin{align*} &v_{2} \left(\frac{D}{\omega } \left(\frac{2}{D_{T}} \right)_{2}
L_{S}(\overline{\psi }_{D_{T}^{3}}, 1) \right) = v_{2}
\left(\frac{D}{\omega } \left(\frac{2}{D_{T}} \right)_{2}
L(\overline{\psi }_{D_{T}^{3}}, 1) \prod _{\pi _{k} \mid \widehat{D}
_{T}} \left(1 - \left(\frac{D_{T}} {\pi _{k}} \right) _{2}
\frac{1}{\pi _{k}} \right) \right) \\
&= v_{2} \left(L (\overline{\psi } _{D_{T} ^{3}}, 1)/ \omega \right)
+ \sum _{\pi _{k} \mid \widehat{D} _{T}} v_{2} \left(1 -
\left(\frac{D_{T}}{\pi _{k}} \right) _{2} \frac{1}{\pi _{k}}
\right). \quad \quad \quad \quad (2.20)
\end{align*}
Note that $ 0 < t(T) < n, $ by induction assumption we have $ v_{2}
\left(L(\overline{\psi }_{D_{T}^{3}}, 1)/ \omega  \right) \geq t(T)
- 1. $ Also $ \left(\frac{D_{T}}{\pi _{k}} \right)_{2} = 1 \
\hbox{or} \ -1 $ for each $ \pi_{k} \mid \widehat{D}_{T}. $ So by
(2.20) above, we get
$$ v_{2} \left(\frac{D}{\omega } \left(\frac{2}{D_{T}} \right)_{2}
L_{S}(\overline{\psi }_{D_{T}^{3}}, 1) \right) \geq t(T) - 1 + n -
t(T) = n - 1. $$ Then together with (2.19) of the case $ T =
\emptyset , $ we obtain
\begin{align*}
&v_{2} \left(L(\overline{\psi }_{D^{3}}, 1) / \omega  \right) =
v_{2} \left(\frac{D}{\omega } \left(\frac{2}{D} \right)_{2}
L_{S}(\overline{\psi }_{D^{3}}, 1) \right) \\
&= v_{2} \left( \left( \sum _{T} \frac{D}{\omega }
\left(\frac{2}{D_{T}} \right)_{2} L_{S}(\overline{\psi }_{D_{T}
^{3}}, 1) \right) - \left(\sum _{T \subsetneqq \{1, \cdots , n \}}
\frac{D}{\omega } \left(\frac{2}{D_{T}} \right)_{2}
L_{S}(\overline{\psi }_{D_{T}^{3}}, 1) \right) \right) \\
&\geq n - 1.
\end{align*}
This proves our conclusion by induction, and the proof is completed.
\quad $ \Box $

\vskip 0.3cm

\hspace{-0.8cm} {\bf References }
\baselineskip 6pt
\parskip 0pt
\vskip 0.4cm
\begin{description}

\item[[BSD]]    B.J.Birch and H.P.F.Swinnerton-Dyer,
{\it {Notes on elliptic curves II}},
J. Reine Angew. Math.  218(1965), 79-108.

\item[[CW]] J. Coates and A. Wiles, on the
conjecture of Birch and Swinnerton-Dyer, Invent. Math.
39 (1977), No.3, 223-251.

\item[[GS]]  C.Coldstein and N.Schappacher,
S$\acute{e}$ries d' Eisenstein et fonction
L de courbes ellipliques $\grave{a}$ multiplication
complexe, J. Reine Agew. Math.,
327(1981), 184-218.

\item[[IR]] K.Ireland  and M.Rosen,
A Classical Introduction to Modern Number Theory, GTM 84, New York:
Springer-Verlag, 1990.

\item[[La]]  S. Lang, Elliptic Functions, second edition,
GTM 112, New York: Springer-Verlag, 1987.

\item[[Law]]  D.F.Lawden, Elliptic Functions and
Applications,  Applied Mathematical Sciences Vol.80, New York:
Springer-Verlag, 1989.

\item[[Le]] F. Lemmermeyer, Reciprocity Laws, New York:
Springer-Verlag, 2000.

\item[[Q]] D. Qiu, On $ p-$adic valuations of $ L(1) $ of elliptic
curves with CM by $ \sqrt{-3}, $ Proceedings of the Royal Society of
Edinburgh, 133A (2003), 1389-1407.

\item[[QZ]] D. Qiu, X. Zhang, Elliptic curves with CM by $ \sqrt{-3}
$ and $ 3-$adic valuations of their $ L-$series, manuscripta math.
108 (2002), 385-397.

\item[[Sil 1]] J.H. Silverman, ``The Arithmetic of Elliptic
Curves'', GTM 106, Springer-Verlag, New York, 1986.

\item[[Sil 2]] J.H. Silverman, ``Advanced Topics in
the Arithmetic of  Elliptic Curves'', GTM 151,
Springer-Verlag, 1994.

\item[[St]]  N.M.Stephens, The diophantine equation
$ x^{3}+y^{3}=Dz^{3} $ and the conjectures
of Birch and Swinnerton-Dyer, J. Reine Angew. Math.,
231(1968), 121-162.

\item[[T]] J. Tate, Algorithm for determining the type of a
singular fiber in an elliptic pencil. In: Modular functions of one
variable, IV (Proc. Internat. Summer School, Univ. Antwert, 1972),
33-52. LNM 476, Springer, Berlin, 1975.

\item[[W]]  A. Weil, Elliptic functions according to
Eisenstein and Kronecker, Springer, 1976.

\item[[Z1]] C. Zhao, A criterion for elliptic curves
with lowest 2-power in L(1),
Math. Proc. Cambridge Philos. Soc. 121(1997), 385-400.

\item[[Z2]] C. Zhao, A criterion for elliptic curves
with second lowest 2-power in L(1), Math. Proc. Cambridge Philos.
Soc. 131(2001), 385-404.

\item[[Z3]]   Chunlai Zhao, A criterion for elliptic curves
with lowest 2-power in L(1) (II), Math. Proc. Cambridge Philos. Soc.
134(2003), 407-420.

\end{description}

\end{document}